# Revisiting two strong approximation results of Dudley and Philipp

*This paper is dedicated to the memory of Walter Philipp.*

**Philippe Berthet[1] and David M. Mason[2],***

*Université Rennes 1 and University of Delaware*

**Abstract:** We demonstrate the strength of a coupling derived from a Gaussian approximation of Zaitsev (1987a) by revisiting two strong approximation results for the empirical process of Dudley and Philipp (1983), and using the coupling to derive extended and refined versions of them.

## 1. Introduction

Einmahl and Mason [17] pointed out in their Fact 2.2 that the Strassen–Dudley theorem (see Theorem 11.6.2 in [11]) in combination with a special case of Theorem 1.1 and Example 1.2 of Zaitsev [42] yields the following coupling. Here $|\cdot|_N$, $N \geq 1$, denotes the usual Euclidean norm on $\mathbb{R}^N$.

**Coupling inequality.** *Let $Y_1, \ldots, Y_n$ be independent mean zero random vectors in $\mathbb{R}^N$, $N \geq 1$, such that for some $B > 0$,*

$$|Y_i|_N \leq B, \ i = 1, \ldots, n.$$

*If $(\Omega, \mathcal{T}, \mathbb{P})$ is rich enough then for each $\delta > 0$, one can define independent normally distributed mean zero random vectors $Z_1, \ldots, Z_n$ with $Z_i$ and $Y_i$ having the same variance/covariance matrix for $i = 1, \ldots, n$, such that for universal constants $C_1 > 0$ and $C_2 > 0$,*

$$(1.1) \qquad \mathbb{P}\left\{\left|\sum_{i=1}^{n}(Y_i - Z_i)\right|_N > \delta\right\} \leq C_1 N^2 \exp\left(-\frac{C_2 \delta}{N^2 B}\right).$$

(Actually Einmahl and Mason did not specify the $N^2$ in (1.1) and they applied a less precise result in [43], however their argument is equally valid when based upon [42].) Often in applications, $N$ is allowed to increase with $n$. This result and its variations, when combined with inequalities from empirical and Gaussian processes and from probability on Banach spaces, has recently been shown to be an extremely powerful tool to establish a Gaussian approximation to the uniform empirical process on the $d-$dimensional cube (Rio [34]), strong approximations for the local empirical process (Einmahl and Mason [17]), extreme value results for the

---

[1]IRMAR, Université Rennes 1, Campus de Beaulieu, 35042 Rennes, France, e-mail: philippe.berthet@univ-rennes1.fr

[2]Statistics Program, University of Delaware, 206 Townsend Hall, Newark, DE 19716, USA, e-mail: davidm@udel.edu

*Research partially supported by an NSF Grant.
*AMS 2000 subject classifications:* primary 62E17, 62E20; secondary 60F15.
*Keywords and phrases:* Gaussian approximation, coupling, strong approximation.

Hopfield model (Bovier and Mason [3] and Gentz and Löwe [19]), laws of the iterated logarithm in Banach spaces (Einmahl and Kuelbs [15]), moderate deviations for Banach space valued sums (Einmahl and Kuelbs [16]), and a functional large deviation result for the local empirical process (Mason [26]). In this paper we shall further demonstrate the strength of (1.1) by revisiting two strong approximation results for the empirical process of Dudley and Philipp [14], and use (1.1) to derive extended and refined versions of them.

Dudley and Philipp [14] was a path breaking paper, which introduced a very effective technique for obtaining Gaussian approximations to sums of i.i.d. Banach space valued random variables. The strong approximation results of theirs, which we shall revisit, were derived from a much more general result in their paper. Key to this result was their Lemma 2.12, which is a special case of an extension by Dehling [8] of a Gaussian approximation in the Prokhorov distance to sums of i.i.d. multivariate random vectors due to Yurinskii [41]. In essence, we shall be substituting the application of their Lemma 2.12 by the above coupling inequality (1.1) based upon Zaitsev [42]. We shall also update and streamline the methodology by employing inequalities that were not available to Dudley and Philipp, when they wrote their paper.

### *1.1. The Gaussian approximation and strong approximation problems*

Let us begin by describing the Gaussian approximation problem for the empirical process. For a fixed integer $n \geq 1$ let $X, X_1, \ldots, X_n$ be independent and identically distributed random variables defined on the same probability space $(\Omega, \mathcal{T}, \mathbb{P})$ and taking values in a measurable space $(\mathcal{X}, \mathcal{A})$. Denote by $\mathbb{E}$ the expectation with respect to $\mathbb{P}$ of real valued random variables defined on $(\Omega, \mathcal{T})$ and write $P = \mathbb{P}^X$. Let $\mathcal{M}$ be the set of all measurable real valued functions on $(\mathcal{X}, \mathcal{A})$. In this paper we consider the following two processes indexed by a sufficiently small class $\mathcal{F} \subset \mathcal{M}$. First, define the *P*-empirical process indexed by $\mathcal{F}$ to be

$$(1.2) \qquad \alpha_n(f) = \frac{1}{\sqrt{n}} \sum_{i=1}^n \{f(X_i) - \mathbb{E}f(X)\}, \ f \in \mathcal{F}.$$

Second, define the *P*-Brownian bridge $\mathbb{G}$ indexed by $\mathcal{F}$ to be the mean zero Gaussian process with the same covariance function as $\alpha_n$,

$$(1.3) \quad \langle f, h \rangle = cov(\mathbb{G}(f), \mathbb{G}(h)) = \mathbb{E}\left(f(X)h(X)\right) - \mathbb{E}\left(f(X)\right)\mathbb{E}(h(X)), \ f, g \in \mathcal{F}.$$

Under entropy conditions on $\mathcal{F}$, the Gaussian process $\mathbb{G}$ has a version which is almost surely continuous with respect to the intrinsic semi-metric

$$(1.4) \qquad d_P(f, h) = \sqrt{\mathbb{E}\left(f(X) - h(X)\right)^2}, \ f, g \in \mathcal{F},$$

that is, we include $d_P$-continuity in the definition of $\mathbb{G}$.

Our goal is to show that a version of $X_1, \ldots, X_n$ and $\mathbb{G}$ can be constructed on the same underlying probability space $(\Omega, \mathcal{T}, \mathbb{P})$ in such a way that

$$(1.5) \qquad \|\alpha_n - \mathbb{G}\|_\mathcal{F} = \sup_{f \in \mathcal{F}} |\alpha_n(f) - \mathbb{G}(f)|$$

is very small with high probability, under useful assumptions on $\mathcal{F}$ and $P$. This is what we call the *Gaussian approximation problem*. We shall use our Gaussian approximation results to define on the same probability $(\Omega, \mathcal{T}, \mathbb{P})$ a sequence $X_1, X_2, \ldots$, i.i.d. $X$ and a sequence $\mathbb{G}_1, \mathbb{G}_2, \ldots$, i.i.d. $\mathbb{G}$ so that with high probability,

$$(1.6) \qquad n^{-1/2} \max_{1 \leq m \leq n} \left\| \sqrt{m}\alpha_m - \sum_{i=1}^m \mathbb{G}_i \right\|_{\mathcal{F}}$$

is small. This is what we call the *strong approximation problem*.

## 1.2. Basic assumptions

We shall assume that $\mathcal{F}$ satisfies the following boundedness condition (F.i) and measurability condition (F.ii).

- **(F.i)** *For some $M > 0$, for all $f \in \mathcal{F}$, $\|f\|_{\mathcal{X}} = \sup_{x \in \mathcal{X}} |f(x)| \leq M/2$.*
- **(F.ii)** *The class $\mathcal{F}$ is point-wise measurable, i.e. there exists a countable subclass $\mathcal{F}_\infty$ of $\mathcal{F}$ such that we can find for any function $f \in \mathcal{F}$ a sequence of functions $\{f_m\}$ in $\mathcal{F}_\infty$ for which $\lim_{m \to \infty} f_m(x) = f(x)$ for all $x \in \mathcal{X}$.*

Assumption (F.i) justifies the finiteness of all the integrals that follow as well as the application of the key inequalities. The requirement (F.ii) is imposed to avoid using outer probability measures in our statements – see Example 2.3.4 in [38].

We intend to compute probability bounds for (1.5) holding for any $n$ and some fixed $M$ in (F.i) with ensuing constants independent of $n$.

## 2. Entropy approach based on Zaitsev [42]

We shall require that one of the following two $L_2$-metric entropy conditions (VC) and (BR) holds on the class $\mathcal{F}$. These conditions are commonly used in the context of weak invariance principles and many examples are available – see e.g. van der Vaart and Wellner [38] and Dudley [12]. In this section we shall state our main results. We shall prove them in Section 5.

## 2.1. $L_2$-covering numbers

First we consider polynomially scattered classes $\mathcal{F}$. Let $F$ be an envelope function for the class $\mathcal{F}$, that is, $F$ a measurable function such that $|f(x)| \leq F(x)$ for all $x \in \mathcal{X}$ and $f \in \mathcal{F}$. Given a probability measure $Q$ on $(\mathcal{X}, \mathcal{A})$ endow $\mathcal{M}$ with the semi-metric $d_Q$, where $d_Q^2(f, h) = \int_{\mathcal{X}} (f-h)^2 dQ$. Further, for any $f \in \mathcal{M}$ set $Q(f^2) = d_Q^2(f, 0) = \int_{\mathcal{X}} f^2 dQ$. For any $\varepsilon > 0$ and probability measure $Q$ denote by $N(\varepsilon, \mathcal{F}, d_Q)$ the minimal number of balls $\{f \in \mathcal{M} : d_Q(f, h) < \varepsilon\}$ of $d_Q$-radius $\varepsilon$ and center $h \in \mathcal{M}$ needed to cover $\mathcal{F}$. The uniform $L_2$-covering number is defined to be

$$(2.1) \qquad N_F(\varepsilon, \mathcal{F}) = \sup_Q N\left(\varepsilon \sqrt{Q(F^2)}, \mathcal{F}, d_Q\right),$$

where the supremum is taken over all probability measures $Q$ on $(\mathcal{X}, \mathcal{A})$ for which $0 < Q(F^2) < \infty$. A class of functions $\mathcal{F}$ satisfying the following uniform entropy condition will be called a VC class.

**(VC)** *Assume that for some $c_0 > 0$, $\nu_0 > 0$, and envelope function $F$,*

$$(2.2) \qquad N_F(\varepsilon, \mathcal{F}) \leq c_0 \varepsilon^{-\nu_0}, \ 0 < \varepsilon < 1.$$

The name "VC class" is given to this condition in recognition to Vapnik and Červonenkis [39] who introduced a condition on classes of sets, which implies (VC). In the sequel we shall assume that $F := M/2$ as in (F.i).

**Proposition 1.** *Under* (F.i), (F.ii) *and* (VC) *with $F := M/2$ for each $\lambda > 1$ there exists a $\rho(\lambda) > 0$ such that for each integer $n \geq 1$ one can construct on the same probability space random vectors $X_1, \ldots, X_n$ i.i.d. $X$ and a version of $\mathbb{G}$ such that*

$$(2.3) \qquad \mathbb{P}\left\{\|\alpha_n - \mathbb{G}\|_{\mathcal{F}} > \rho(\lambda) n^{-\tau_1} (\log n)^{\tau_2}\right\} \leq n^{-\lambda},$$

*where $\tau_1 = 1/(2 + 5\nu_0)$ and $\tau_2 = (4 + 5\nu_0)/(4 + 10\nu_0)$.*

Proposition 1 leads to the following strong approximation result. It is an indexed by functions generalization of an indexed by sets result given in Theorem 7.4 of Dudley and Philipp [14].

**Theorem 1.** *Under the assumptions and notation of Proposition 1 for all $1/(2\tau_1) < \alpha < 1/\tau_1$ and $\gamma > 0$ there exist a $\rho(\alpha, \gamma) > 0$, a sequence of i.i.d. $X_1, X_2, \ldots$, and a sequence of independent copies $\mathbb{G}_1, \mathbb{G}_2, \ldots$, of $\mathbb{G}$ sitting on the same probability space such that*

$$(2.4) \qquad \mathbb{P}\left\{\max_{1 \leq m \leq n} \left\|\sqrt{m}\alpha_m - \sum_{i=1}^{m} \mathbb{G}_i\right\|_{\mathcal{F}} > C\rho(\alpha, \gamma) n^{1/2 - \tau(\alpha)} (\log n)^{\tau_2}\right\} \leq n^{-\gamma}$$

*and*

$$(2.5) \qquad \max_{1 \leq m \leq n} \left\|\sqrt{m}\alpha_m - \sum_{i=1}^{m} \mathbb{G}_i\right\|_{\mathcal{F}} = O\left(n^{1/2 - \tau(\alpha)} (\log n)^{\tau_2}\right), \ a.s.,$$

*where $\tau(\alpha) = (\alpha\tau_1 - 1/2)/(1 + \alpha) > 0$.*

### 2.2. Bracketing numbers

A second way to measure the size of the class $\mathcal{F}$ is to use $L_2(P)$-brackets instead of $L_2(Q)$-balls. Let $l \in \mathcal{M}$ and $u \in \mathcal{M}$ be such that $l \leq u$ and $d_P(l, u) < \varepsilon$. The pair of functions $l, u$ form an $\varepsilon$-bracket $[l, u]$ consisting of all the functions $f \in \mathcal{F}$ such that $l \leq f \leq u$. Let $N_{[\,]}(\varepsilon, \mathcal{F}, d_P)$ be the minimum number of $\varepsilon$-brackets needed to cover $\mathcal{F}$. Notice that trivially we have $N(\varepsilon, \mathcal{F}, d_P) \leq N_{[\,]}(\varepsilon/2, \mathcal{F}, d_P)$.

**(BR)** *Assume that for some $b_0 > 0$ and $0 < r_0 < 1$,*

$$(2.6) \qquad \log N_{[\,]}(\varepsilon, \mathcal{F}, d_P) \leq b_0^2 \varepsilon^{-2r_0}, \ 0 < \varepsilon < 1.$$

We derive the following rate of Gaussian approximation assuming an exponentially scattered index class $\mathcal{F}$, meaning that (2.6) holds. Note that we get a slower rate in Proposition 2 than that given Proposition 1.

**Proposition 2.** *Under* (F.i), (F.ii) *and* (BR) *for each* $\lambda > 1$ *there exists a* $\rho(\lambda) > 0$ *such that for each integer* $n \geq 1$ *one can construct on the same probability space random vectors* $X_1, \ldots, X_n$ *i.i.d.* $X$ *and a version of* $\mathbb{G}$ *such that*

$$\mathbb{P}\left\{\|\alpha_n - \mathbb{G}\|_\mathcal{F} > \rho(\lambda)(\log n)^{-\kappa}\right\} \leq n^{-\lambda}, \tag{2.7}$$

*where* $\kappa = (1 - r_0)/2r_0$.

Proposition 2 leads to the following indexed by functions generalization of an indexed by sets result given in Theorem 7.1 of Dudley and Philipp [14].

**Theorem 2.** *Under the assumptions and notation of Proposition 2, with* $\kappa < 1/2$ *($1/2 < r_0 < 1$), for every* $H > 0$ *there exist* $\rho(\tau, H) > 0$ *and a sequence of i.i.d.* $X_1, X_2, \ldots$, *and a sequence of independent copies* $\mathbb{G}_1, \mathbb{G}_2, \ldots$, *of* $\mathbb{G}$ *sitting on the same probability space such that*

$$\mathbb{P}\left\{\max_{1 \leq m \leq n}\left\|\sqrt{m}\alpha_m - \sum_{i=1}^m \mathbb{G}_i\right\|_\mathcal{F} > \sqrt{n}\rho(\tau, H)(\log n)^{-\tau}\right\} \leq (\log n)^{-H} \tag{2.8}$$

*and*

$$\max_{1 \leq m \leq n}\left\|\sqrt{m}\alpha_m - \sum_{i=1}^m \mathbb{G}_i\right\|_\mathcal{F} = O\left(\sqrt{n}(\log n)^{-\tau}\right), \ a.s., \tag{2.9}$$

*where* $\tau = \kappa(1/2 - \kappa)/(1 - \kappa)$.

## 3. Comments on the approach based on KMT

Given $\mathcal{F}$, the rates obtained in Proposition 1 and Theorem 1 are universal in $P$. If one specializes to particular $P$, the rates in Propositions 1 and 2 and Theorem 1 and 2 are far from being optimal. In such situations one can get better and even unimprovable rates by replacing the use of Zaitsev [42] by the Komlós, Major and Tusnády [KMT] [22] Brownian bridge approximation to the uniform empirical process or one based on the same dyadic scheme. (More details about this approximation are provided in [4, 13, 25, 27, 28].) This is especially the case when the underlying probability measure $P$ is smooth. To see how this works in the empirical process indexed by functions setup refer to Koltchinskii [21] and Rio [33] and in the indexed by smooth sets situation turn to Révesz [32] and Massart [29]. One can also use the KMT–type bivariate Brownian bridge approximation to the bivariate uniform empirical process as a basis for further approximation. For a brief outline of this approximation consult Tusnády [36] and for detailed presentations refer to Castelle [5] and Castelle and Laurent-Bonvalot [6].

## 4. Tools needed in proofs

For convenience we shall collect here the basic tools we shall need in our proofs.

### *4.1. Inequalities for empirical processes*

On a rich enough probability space $(\Omega, \mathcal{T}, \mathbb{P})$, let $X, X_1, X_2, \ldots, X_n$ be i.i.d. random variables with law $P = \mathbb{P}^X$ and $\epsilon_1, \epsilon_2, \ldots, \epsilon_n$ be i.i.d. Rademacher random variables independent of $X_1, \ldots, X_n$. By a Rademacher random variable $\epsilon_1$, we mean that $\mathbb{P}(\epsilon_1 = 1) = \mathbb{P}(\epsilon_1 = -1) = 1/2$. Consider a point-wise measurable class $\mathcal{G}$ of bounded measurable real valued functions on $(\mathcal{X}, \mathcal{A})$.

The following exponential inequality is due to Talagrand [35].

**Talagrand's inequality.** *If $\mathcal{G}$ satisfies* (F.i) *and* (F.ii) *then for all $n \geq 1$ and $t > 0$ we have, for suitable finite constants $A > 0$ and $A_1 > 0$,*

$$\begin{aligned}
&\mathbb{P}\left\{\|\alpha_n\|_{\mathcal{G}} > A\left(\mathbb{E}\left(\left\|\frac{1}{\sqrt{n}}\sum_{i=1}^n \epsilon_i g(X_i)\right\|_{\mathcal{G}}\right) + t\right)\right\} \\
&\qquad \leq 2\exp\left(-\frac{A_1 t^2}{\sigma_{\mathcal{G}}^2}\right) + 2\exp\left(-\frac{A_1 t \sqrt{n}}{M}\right),
\end{aligned} \quad (4.1)$$

*where $\sigma_{\mathcal{G}}^2 := \sup_{g \in \mathcal{G}} Var(g(X))$.*

Moreover the constants $A$ and $A_1$ are independent of $\mathcal{G}$ and $M$. Next we state two upper bounds for the above expectation of the supremum of the symmetrized empirical process.

We shall require two moment bounds. The first is due to Einmahl and Mason [18] – for a similar bound refer to Giné and Guillou [20].

**Moment inequality for (VC).** *Let $\mathcal{G}$ satisfy* (F.i) *and* (F.ii) *with envelope function $G$ and be such that for some positive constants $\beta, v, c > 1$ and $\sigma \leq 1/(8c)$ the following four conditions hold,*

$$\mathbb{E}(G^2(X)) \leq \beta^2; \ N_G(\varepsilon, \mathcal{G}) \leq c\varepsilon^{-v}, \ 0 < \varepsilon < 1; \ \sup_{g \in \mathcal{G}} \mathbb{E}(g^2(X)) \leq \sigma^2;$$

*and*

$$\sup_{g \in \mathcal{G}} \|g\|_{\mathcal{X}} \leq \frac{\sqrt{n\sigma^2/\log(\beta \vee 1/\sigma)}}{2\sqrt{v+1}}.$$

*Then we have for a universal constant $A_2$ not depending on $\beta$,*

$$\mathbb{E}\left(\left\|\frac{1}{\sqrt{n}}\sum_{i=1}^n \epsilon_i g(X_i)\right\|_{\mathcal{G}}\right) \leq A_2\sqrt{v\sigma^2 \log(\beta \vee 1/\sigma)}. \quad (4.2)$$

Next we state a moment inequality under (BR). For any $0 < \sigma < 1$, set

$$J(\sigma, \mathcal{G}) = \int_{[0,\sigma]} \sqrt{\log N_{[\,]}(s, \mathcal{G}, d_P)}\, ds \quad (4.3)$$

and

$$a(\sigma, \mathcal{G}) = \frac{\sigma}{\sqrt{\log N_{[\,]}(\sigma, \mathcal{G}, d_P)}}. \quad (4.4)$$

The second moment bound follows from Lemma 19.34 in [37] and a standard symmetrization inequality, and is reformulated by using (4.3).

**Moment inequality for (BR).** *Let $\mathcal{G}$ satisfy* (F.i) *and* (F.ii) *with envelope $G$ and be such that $\sup_{g \in \mathcal{G}} \mathbb{E}(g^2(X)) < \sigma^2 < 1$. We have, for a universal constant $A_3$,*

$$\mathbb{E}\left(\left\|\frac{1}{\sqrt{n}}\sum_{i=1}^n \epsilon_i g(X_i)\right\|_{\mathcal{G}}\right) \leq A_3\left(J(\sigma, \mathcal{G}) + \sqrt{n}\,\mathbb{P}\{G(X) > \sqrt{n}\,a(\sigma, \mathcal{G})\}\right). \quad (4.5)$$

### 4.2. Inequalities for Gaussian processes

Let $\mathbb{Z}$ be a separable mean zero Gaussian process on a probability space $(\Omega, \mathcal{T}, \mathbb{P})$ indexed by a set $T$. Define the intrinsic semi–metric $\rho$ on $T$ by

$$(4.6) \qquad \rho(s,t) = \sqrt{\mathbb{E}(\mathbb{Z}_t - \mathbb{Z}_s)^2}.$$

For each $\varepsilon > 0$ let $N(\varepsilon, T, \rho)$ denote the minimal number of $\rho$-balls of radius $\varepsilon$ needed to cover $T$. Write $\|\mathbb{Z}\|_T = \sup_{t \in T} |\mathbb{Z}_t|$ and $\sigma_T^2(\mathbb{Z}) = \sup_{t \in T} \mathbb{E}(\mathbb{Z}_t^2)$. The following large deviation probability estimate for $\|\mathbb{Z}\|_T$ is due to Borell [2]. (Also see Proposition A.2.1 in [38].)

**Borell's inequality.** *For all $t > 0$,*

$$(4.7) \qquad \mathbb{P}\{|\|\mathbb{Z}\|_T - \mathbb{E}(\|\mathbb{Z}\|_T)| > t\} \leq 2\exp\left(-\frac{t^2}{2\sigma_T^2(\mathbb{Z})}\right).$$

According to Dudley [9], the entropy condition

$$(4.8) \qquad \int_{[0,1]} \sqrt{\log N(\varepsilon, T, \rho)}\, d\varepsilon < \infty$$

ensures the existence of a separable, bounded, $d_P$-uniformly continuous modification of $\mathbb{Z}$. Moreover the above Dudley integral (4.8) controls the modulus of continuity of $\mathbb{Z}$ (see Dudley [10]) as well as its expectation (see Marcus and Pisier [24], p. 25, Ledoux and Talagrand [23], p. 300, de la Peña and Giné [7], Cor. 5.1.6, and Dudley [12]). The following inequality is part of Corollary 2.2.8 in van der Vaart and Wellner [38].

**Gaussian moment inequality.** *For some universal constant $A_4 > 0$ and all $\sigma > 0$ we have*

$$(4.9) \qquad \mathbb{E}\left(\sup_{\rho(s,t)<\sigma} |\mathbb{Z}_t - \mathbb{Z}_s|\right) \leq A_4 \int_{[0,\sigma]} \sqrt{\log N(\varepsilon, T, \rho)}\, d\varepsilon.$$

We shall be applying these inequalities to the Gaussian process $\mathbb{Z} = \mathbb{G}$ defined in introduction, so that $T = \mathcal{F}$ and $\rho = d_P$.

### 4.3. A maximal inequality

The following version of a maximal inequality due to Montgomery–Smith [30] (see also Theorem 1.1.5 in [7]) will come in handy.

**A maximal inequality.** *Let $X_1, \ldots, X_n$, $n \geq 1$, be i.i.d. random variables taking values in a separable Banach space. Then for all $t > 0$,*

$$(4.10) \qquad \mathbb{P}\left\{\max_{1 \leq m \leq n} \left\|\sum_{i=1}^m X_i\right\| > t\right\} \leq 9\mathbb{P}\left\{\left\|\sum_{i=1}^n X_i\right\| > \frac{t}{30}\right\}.$$

## 5. Proofs of main results

### 5.1. Description of construction of $(\alpha_n, \mathbb{G})$

Under (F.i), (F.ii) and either (VC) or (BR) for any $\varepsilon > 0$ we can choose a grid

$$\mathcal{H}(\varepsilon) = \{h_k : 1 \leq k \leq N(\varepsilon)\}$$

of measurable functions on $(\mathcal{X}, \mathcal{A})$ such that each $f \in \mathcal{F}$ is in a ball $\{f \in \mathcal{M} : d_P(h_k, f) < \varepsilon\}$ around some $h_k$, $1 \leq k \leq N(\varepsilon)$. The choice

(5.1) $$N(\varepsilon) \leq N(\varepsilon/2, \mathcal{F}, d_P)$$

permits us to select $h_k \in \mathcal{F}$. Set

$$\mathcal{F}(\varepsilon) = \{(f, f') \in \mathcal{F}^2 : d_P(f, f') < \varepsilon\}.$$

Fix $n \geq 1$. Let $X, X_1, \ldots, X_n$ be independent with common law $P = \mathbb{P}^X$ and $\epsilon_1, \ldots, \epsilon_n$ be independent Rademacher random variables mutually independent of $X_1, \ldots, X_n$. Write for $\varepsilon > 0$,

$$\mu_n(\varepsilon) = \mathbb{E}\left\{\sup_{(f,f') \in \mathcal{F}(\varepsilon)} \left|\frac{1}{\sqrt{n}} \sum_{i=1}^n \epsilon_i (f - f')(X_i)\right|\right\}$$

and

$$\mu(\varepsilon) = \mathbb{E}\left\{\sup_{(f,f') \in \mathcal{F}(\varepsilon)} |\mathbb{G}(f) - \mathbb{G}(f')|\right\}.$$

Given $\varepsilon > 0$ and $n \geq 1$, our aim is to construct a probability space $(\Omega, \mathcal{T}, \mathbb{P})$ on which sit $X_1, \ldots, X_n$ and a version of the Gaussian process $\mathbb{G}$ indexed by $\mathcal{F}$ such that for $\mathcal{H}(\varepsilon)$ and $\mathcal{F}(\varepsilon)$ defined as above and for all $A > 0$, $\delta > 0$ and $t > 0$,

$$\mathbb{P}\{\|\alpha_n - \mathbb{G}\|_{\mathcal{F}} > A\mu_n(\varepsilon) + \mu(\varepsilon) + \delta + (A+1)t\}$$

$$\leq \mathbb{P}\left\{\max_{h \in \mathcal{H}(\varepsilon)} |\alpha_n(h) - \mathbb{G}(h)| > \delta\right\}$$

(5.2) $$+ \mathbb{P}\left\{\sup_{(f,f') \in \mathcal{F}(\varepsilon)} |\alpha_n(f) - \alpha_n(f')| > A\mu_n(\varepsilon) + At\right\}$$

$$+ \mathbb{P}\left\{\sup_{(f,f') \in \mathcal{F}(\varepsilon)} |\mathbb{G}(f) - \mathbb{G}(f')| > t + \mu(\varepsilon)\right\}$$

$$=: P_n(\delta) + Q_n(t, \varepsilon) + Q(t, \varepsilon),$$

with all these probabilities simultaneously small for suitably chosen $A > 0$, $\delta > 0$ and $t > 0$. Consider the $n$ i.i.d. mean zero random vectors in $\mathbb{R}^{N(\varepsilon)}$,

$$Y_i := \frac{1}{\sqrt{n}} \left(h_1(X_i) - \mathbb{E}(h_1(X)), \ldots, h_{N(\varepsilon)}(X_i) - \mathbb{E}(h_{N(\varepsilon)}(X))\right), \ 1 \leq i \leq n.$$

First note that by $h_k \in \mathcal{F}$ and (F.i), we have

$$|Y_i|_{N(\varepsilon)} \leq M\sqrt{\frac{N(\varepsilon)}{n}}, \ 1 \leq i \leq n.$$

Therefore by the coupling inequality (1.1) we can define $Y_1, \ldots, Y_n$ i.i.d.

$$Y := \left(Y^1, \ldots, Y^{N(\varepsilon)}\right)$$

and $Z_1, \ldots, Z_n$ i.i.d.

$$Z := \left(Z^1, \ldots, Z^{N(\varepsilon)}\right)$$

mean zero Gaussian vectors on the same probability space such that

$$(5.3) \quad P_n(\delta) \leq \mathbb{P}\left\{\left|\sum_{i=1}^n (Y_i - Z_i)\right|_{N(\varepsilon)} > \delta\right\} \leq C_1 N(\varepsilon)^2 \exp\left(-\frac{C_2 \sqrt{n}\, \delta}{(N(\varepsilon))^{5/2} M}\right),$$

where $cov(Z^l, Z^k) = cov(Y^l, Y^k) = \langle h_l, h_k \rangle$. Moreover by Lemma A1 of Berkes and Philipp [1] (also see Vorob'ev [40]) this space can be extended to include a $P$-Brownian bridge $\mathbb{G}$ indexed by $\mathcal{F}$ such that

$$\mathbb{G}(h_k) = n^{-1/2} \sum_{i=1}^n Z_i^k.$$

The $P_n(\delta)$ in (5.2) is defined through this $\mathbb{G}$. Notice that the probability space on which $Y_1, \ldots, Y_n, Z_1, \ldots, Z_n$ and $\mathbb{G}$ sit depends on $n \geq 1$ and the choice of $\varepsilon > 0$ and $\delta > 0$.

Observe that the class

$$\mathcal{G}(\varepsilon) = \{f - f' : (f, f') \in \mathcal{F}(\varepsilon)\}$$

satisfies (F.i) with $M/2$ replaced by $M$, (F.ii) and

$$\sigma_{\mathcal{G}(\varepsilon)}^2 = \sup_{(f,f')\in\mathcal{F}(\varepsilon)} Var(f(X) - f'(X)) \leq \sup_{(f,f')\in\mathcal{F}(\varepsilon)} d_P^2(f, f') \leq \varepsilon^2.$$

Thus with $A > 0$ as in (4.1) we get by applying Talagrand's inequality,

$$(5.4) \quad \begin{aligned} Q_n(t, \varepsilon) &= \mathbb{P}\left\{\|\alpha_n\|_{\mathcal{G}(\varepsilon)} > A(\mu_n(\varepsilon) + t)\right\} \\ &\leq 2\exp\left(-\frac{A_1 t^2}{\varepsilon^2}\right) + 2\exp\left(-\frac{A_1 \sqrt{n}\, t}{M}\right). \end{aligned}$$

Next, consider the separable centered Gaussian process $\mathbb{Z}_{(f,f')} = \mathbb{G}(f) - \mathbb{G}(f')$ indexed by $T = \mathcal{F}(\varepsilon)$. We have

$$\begin{aligned} \sigma_T^2(\mathbb{Z}) &= \sup_{(f,f')\in\mathcal{F}(\varepsilon)} \mathbb{E}\left((\mathbb{G}(f) - \mathbb{G}(f'))^2\right) = \sup_{(f,f')\in\mathcal{F}(\varepsilon)} Var(f(X) - f'(X)) \\ &\leq \sup_{(f,f')\in\mathcal{F}(\varepsilon)} d_P^2(f, f') \leq \varepsilon^2. \end{aligned}$$

Borell's inequality (4.7) now gives

$$(5.5) \quad Q(t, \varepsilon) = \mathbb{P}\left\{\sup_{(f,f')\in\mathcal{F}(\varepsilon)} |\mathbb{G}(f) - \mathbb{G}(f')| > t + \mu(\varepsilon)\right\} \leq 2\exp\left(-\frac{t^2}{2\varepsilon^2}\right).$$

Putting (5.3), (5.4) and (5.5) together we obtain, for some positive constants $A$, $A_1$ and $A_5$ with $A_5 \leq 1/2$,

(5.6)
$$\mathbb{P}\{\|\alpha_n - \mathbb{G}\|_{\mathcal{F}} > A\mu_n(\varepsilon) + \mu(\varepsilon) + \delta + (A+1)t\}$$
$$\leq C_1 N(\varepsilon)^2 \exp\left(-\frac{C_2\sqrt{n}\,\delta}{(N(\varepsilon))^{5/2} M}\right)$$
$$+ 2\exp\left(-\frac{A_1\sqrt{n}\,t}{M}\right) + 4\exp\left(-\frac{A_5 t^2}{\varepsilon^2}\right).$$

*Proof of Proposition 1.* Let us assume that (VC) holds with $F := M/2$, so that for some $c_0 > 0$ and $\nu_0 > 0$, with $c_1 = c_0(2\sqrt{PF^2})^{\nu_0} = c_0 M^{\nu_0}$,

$$N(\varepsilon) \leq N(\varepsilon/2, \mathcal{F}, d_P) \leq c_1 \varepsilon^{-\nu_0},\ 0 < \varepsilon < 1.$$

Notice that both

$$N(\varepsilon, \mathcal{G}(\varepsilon), d_P) \leq (N(\varepsilon/2, \mathcal{F}, d_P))^2 \leq c_1^2 \varepsilon^{-2\nu_0}$$

and

$$N(\varepsilon, \mathcal{F}(\varepsilon), d_P) \leq (N(\varepsilon/2, \mathcal{F}, d_P))^2 \leq c_1^2 \varepsilon^{-2\nu_0}.$$

Therefore we can apply the moment bound assuming (VC) given in (4.2) taken with $\mathcal{G} = \mathcal{G}(\varepsilon)$, $G := M$, $\upsilon = 2\nu_0$ and $\beta = M$, to get for any $0 < \varepsilon < 1/e$ and $n \geq 1$ so that

(5.7)
$$\frac{\sqrt{n}\varepsilon}{2\sqrt{1+2\nu_0}\sqrt{\log(M \vee 1/\varepsilon)}} > M$$

the bound

$$\mu_n(\varepsilon) \leq A_2 \varepsilon \sqrt{2\nu_0 \log(M \vee 1/\varepsilon)}.$$

Whereas, by the Gaussian moment bound (4.9), we have for all $0 < \varepsilon < 1/e$,

$$\mu(\varepsilon) \leq A_4 \sqrt{2\nu_0} \int_{[0,\varepsilon]} \sqrt{\log(1/x)}\,dx.$$

Hence, for some $D > 0$ it holds for all $0 < \varepsilon < 1/e$ and $n \geq 1$ so that (5.7) holds,

(5.8)
$$A\mu_n(\varepsilon) + \mu(\varepsilon) \leq D\varepsilon\sqrt{\log(1/\varepsilon)}.$$

Therefore, in view of (5.8) and (5.6) it is natural to define for suitably large positive $\gamma_1$ and $\gamma_2$,

$$\delta = \gamma_1 \varepsilon \sqrt{\log(1/\varepsilon)} \text{ and } t = \gamma_2 \varepsilon \sqrt{\log(1/\varepsilon)}.$$

We now have for all $0 < \varepsilon < 1/e$ and $n \geq 1$ so that (5.7) is satisfied on a suitable probability space depending on $n \geq 1$, $\varepsilon$ and $\delta$ so that (5.6) holds,

$$\mathbb{P}\left\{\|\alpha_n - \mathbb{G}\|_{\mathcal{F}} > (D + \gamma_1 + (1+A)\gamma_2)\varepsilon\sqrt{\log(1/\varepsilon)}\right\}$$
$$\leq \frac{C_1 c_1^2}{\varepsilon^{2\nu_0}} \exp\left(-\frac{\gamma_1 C_2 \sqrt{n}}{c_1^{5/2} M}\varepsilon^{1+5\nu_0/2}\sqrt{\log(1/\varepsilon)}\right)$$
$$+ 2\exp\left(-\frac{A_1 \gamma_2 \sqrt{n}}{M}\varepsilon\sqrt{\log(1/\varepsilon)}\right) + 4\exp\left(-A_5 \gamma_2^2 \log(1/\varepsilon)\right).$$

By taking $\varepsilon = ((\log n)/n)^{1/(2+5\nu_0)}$, which satisfies (5.7) for all large enough $n$, we readily obtain from these last bounds that for every $\lambda > 1$ there exist $D > 0$, $\gamma_1 > 0$ and $\gamma_2 > 0$ such that for all $n \geq 1$, $\alpha_n$ and $\mathbb{G}$ can be defined on the same probability space so that

$$\mathbb{P}\left\{\|\alpha_n - \mathbb{G}\|_{\mathcal{F}} > (D + \gamma_1 + (1+A)\gamma_2)\left(\frac{\log n}{n}\right)^{1/(2+5\nu_0)}\sqrt{\frac{\log n}{2+5\nu_0}}\right\} \leq n^{-\lambda}.$$

It is clear now that there exists a $\rho(\lambda) > 0$ such that (2.3) holds. This completes the proof of Proposition 1. □

*Proof of Proposition 2.* Under (BR) as defined in (2.6) we have, for some $0 < r_0 < 1$ and $b_0 > 0$,

$$N(\varepsilon) \leq N(\varepsilon/2, \mathcal{F}, d_P) \leq N_{[\,]}(\varepsilon/2, \mathcal{F}, d_P) \leq \exp\left(\frac{2^{2r_0}b_0^2}{\varepsilon^{2r_0}}\right), \quad 0 < \varepsilon < 1,$$

and as above both

$$N(\varepsilon, \mathcal{G}(\varepsilon), d_P) \leq N_{[\,]}(\varepsilon, \mathcal{G}(\varepsilon), d_P) \leq \left(N_{[\,]}(\varepsilon/2, \mathcal{F}, d_P)\right)^2 \leq \exp\left(2\frac{2^{2r_0}b_0^2}{\varepsilon^{2r_0}}\right)$$

and

$$N(\varepsilon, \mathcal{F}(\varepsilon), d_P) \leq N_{[\,]}(\varepsilon, \mathcal{F}(\varepsilon), d_P) \leq \left(N_{[\,]}(\varepsilon/2, \mathcal{F}, d_P)\right)^2 \leq \exp\left(2\frac{2^{2r_0}b_0^2}{\varepsilon^{2r_0}}\right).$$

Setting $\sigma = \varepsilon$ in (4.3) and (4.4) we get

$$J(\varepsilon, \mathcal{G}(\varepsilon)) \leq \sqrt{2}b_0 \int_{[0,\varepsilon]} \frac{ds}{s^{r_0}} \leq \frac{\sqrt{2}b_0}{1-r_0}\varepsilon^{1-r_0}$$

and

$$a(\varepsilon, \mathcal{G}(\varepsilon)) = \frac{\varepsilon}{\sqrt{\log N_{[\,]}(\varepsilon, \mathcal{G}(\varepsilon), d_P)}} \geq \frac{\varepsilon^{1+r_0}}{\sqrt{2}b_0}.$$

Hence by the moment bound assuming (BR) given in (4.5) taken with $G(X) = M$,

$$\mu_n(\varepsilon) \leq A_3\left(\frac{\sqrt{2}b_0}{1-r_0}\varepsilon^{1-r_0} + \sqrt{n}\,\mathbb{I}_{\left\{M > \frac{\sqrt{n}\varepsilon^{1+r_0}}{\sqrt{2}b_0}\right\}}\right)$$

and, since in the same way we have

$$J(\varepsilon, \mathcal{F}(\varepsilon)) \leq \frac{\sqrt{2}b_0}{1-r_0}\varepsilon^{1-r_0} \text{ and } a(\varepsilon, \mathcal{F}(\varepsilon)) \geq \frac{\varepsilon^{1+r_0}}{\sqrt{2}b_0},$$

we get by the Gaussian moment inequality,

$$\mu(\varepsilon) \leq \frac{A_4\sqrt{2}b_0}{1-r_0}\varepsilon^{1-r_0}.$$

As a consequence, for some $D > 0$ and

$$\varepsilon > \frac{(DM)^{1/(1+r_0)}}{n^{1/(2+2r_0)}}$$

it follows that
$$A\mu_n(\varepsilon) + \mu(\varepsilon) \le D\varepsilon^{1-r_0}.$$

Thus it is natural to take in (5.6) for some $\gamma_1 > 0$ and $\gamma_2 > 0$ large enough,
$$\delta = \gamma_1 \varepsilon^{1-r_0} \text{ and } t = \gamma_2 \varepsilon^{1-r_0},$$

which gives with $\rho = D + \gamma_1 + (A+1)\gamma_2$,

$$\mathbb{P}\left\{\|\alpha_n - \mathbb{G}\|_{\mathcal{F}} > \rho\varepsilon^{1-r_0}\right\}$$
$$\le C_1 \exp\left(\frac{2^{2r_0+1}b_0^2}{\varepsilon^{2r_0}} - \frac{\gamma_1 C_2 \sqrt{n}}{M}\varepsilon^{1-r_0}\exp\left(-\frac{5\left(2^{2r_0}b_0^2\right)}{2\varepsilon^{2r_0}}\right)\right)$$
$$+ 2\exp\left(-\frac{A_1\gamma_2\sqrt{n}}{M}\varepsilon^{1-r_0}\right) + 4\exp\left(-\frac{A_5\gamma_2^2}{\varepsilon^{2r_0}}\right).$$

We choose
$$\varepsilon = \left(\frac{10b_0^2 2^{2r_0}}{\log n}\right)^{1/(2r_0)},$$
which makes
$$\exp\left(-\frac{5\left(2^{2r_0}b_0^2\right)}{2\varepsilon^{2r_0}}\right) = n^{-1/4}.$$

Given any $\lambda > 0$ we clearly see now from this last probability bound that for $\rho(\lambda) > 0$ made large enough by increasing $\gamma_1$ and $\gamma_2$ we get for all $n \ge 1$,
$$\mathbb{P}\left\{\|\alpha_n - \mathbb{G}\|_{\mathcal{F}} > \rho(\lambda)(\log n)^{-(1-r_0)/2r_0}\right\} \le n^{-\lambda}.$$

The proof of Proposition 2 now follows the same lines as that of Proposition 1. □

### 5.2. Proofs of strong approximations

Notice that the conditions on $\mathcal{F}$ in Propositions 1 and 2 imply that there exists a constant $B$ such that
$$\sup_{n\ge 1}\mathbb{E}\left(\left\|\frac{1}{\sqrt{n}}\sum_{i=1}^n \epsilon_i f(X_i)\right\|_{\mathcal{F}}\right) \le B \text{ and } \mathbb{E}(\|\mathbb{G}\|_{\mathcal{F}}) \le B.$$

Therefore by Talagrand's inequality (4.1) and the Montgomery–Smith inequality (4.10) for all $n \ge 1$ and $t > 0$ we have, for suitable finite constants $C > 0$ and $C_1 > 0$,

(5.9)
$$\mathbb{P}\left\{\max_{1\le m\le n}\sqrt{m}\|\alpha_m\|_{\mathcal{F}} > C\sqrt{n}(B+t)\right\}$$
$$\le 18\exp\left(-\frac{C_1 t^2}{\sigma_{\mathcal{F}}^2}\right) + 18\exp\left(-\frac{C_1 t\sqrt{n}}{M}\right),$$

where $\sigma_{\mathcal{F}}^2 := \sup_{f\in\mathcal{F}} Var(f(X))$. Furthermore, by Borell's inequality (4.7), the Montgomery–Smith inequality (4.10) and the fact that $n^{-1/2}\sum_{i=1}^n \mathbb{G}_i =_d \mathbb{G}$, for i.i.d. $\mathbb{G}_i$, we get for all $n \ge 1$ and $t > 0$ that for a suitable finite constant $D > 0$,

(5.10)    $$\mathbb{P}\left\{\max_{1\le m\le n}\left\|\sum_{i=1}^m \mathbb{G}_i\right\|_{\mathcal{F}} > D\sqrt{n}(B+t)\right\} \le 18\exp\left(-\frac{t^2}{2\sigma_{\mathcal{F}}^2}\right).$$

*Proof of Theorem 1.* Choose any $\gamma > 0$. We shall modify the scheme described on pages 236–238 of Philipp [31] to construct a probability space on which (2.4) and (2.5) hold. Let $n_0 = 1$ and for each $k \geq 1$ set $n_k = [k^\alpha]$, where $[x]$ denotes the integer part of $x$ and $\alpha$ is chosen so that

(5.11) $$1/2 < \tau_1 \alpha < 1.$$

Notice that $\tau_1 < 1/2$ in Proposition 1 and thus $\alpha > 1$.

Applying Proposition 1, we see that for each $\lambda > 1$ there exists a $\rho = \rho(\lambda) > 0$ such that one can construct a sequence of independent pairs $(\alpha_{n_k}^{(k)}, \mathbb{G}^{(k)})_{k \geq 1}$ sitting on the same probability space satisfying for all $k \geq 1$,

(5.12) $$\mathbb{P}\left\{ \left\| \alpha_{n_k}^{(k)} - \mathbb{G}^{(k)} \right\|_{\mathcal{F}} > \rho n_k^{-\tau_1} (\log n_k)^{\tau_2} \right\} \leq n_k^{-\lambda}.$$

Set for $k \geq 1$
$$t_k = \sum_{j < k} n_j \sim \frac{1}{1+\alpha} k^{\alpha+1}.$$

Using Lemma A1 of Berkes and Philipp [1] we can assume that each $\alpha_{n_k}^{(k)}$ is formed from $X_{t_k+1}, \ldots, X_{t_{k+1}}$ i.i.d. $X$ and that each $\mathbb{G}^{(k)}$ is formed as

$$\mathbb{G}^{(k)} = \frac{1}{\sqrt{n_k}} \sum_{t_k < j \leq t_{k+1}} \mathbb{G}_j,$$

where $\mathbb{G}_{t_k+1}, \ldots, \mathbb{G}_{t_{k+1}}$ are i.i.d. $\mathbb{G}$. Moreover we can do this in such a way that $X_1, X_2 \ldots,$ are i.i.d. $X$ and $\mathbb{G}_1, \mathbb{G}_2, \ldots,$ are i.i.d. $\mathbb{G}$. For any integer $N \geq 2$ set $N(\beta) = [N^\beta]$, where $\beta = \alpha/(1+\alpha)$. Define

$$s(N) = \sum_{k=N(\beta)}^{N} n_k^{1/2-\tau_1} (\log n_k)^{\tau_2}.$$

Now for some constants $c_1 > 0$ and $c > 0$,

(5.13) $$s(N) \sim c_1 N^{(1+\alpha)/2 - (\alpha \tau_1 - 1/2)} (\log N)^{\tau_2} \sim c (t_N)^{1/2 - \tau(\alpha)} (\log t_N)^{\tau_2},$$

where $\tau(\alpha) = (\alpha \tau_1 - 1/2)/(1+\alpha) > 0$, by (5.11).

We have

$$\mathbb{P}\left\{ \max_{1 \leq m \leq t_N} \left\| \sum_{j=1}^{m} [f(X_j) - \mathbb{E}f(X) - \mathbb{G}_j(f)] \right\|_{\mathcal{F}} > \rho s(N) \right\}$$

$$\leq \mathbb{P}\left\{ \max_{1 \leq m \leq t_{N(\beta)}} \left\| \sum_{j=1}^{m} [f(X_j) - \mathbb{E}f(X)] \right\|_{\mathcal{F}} > \frac{\rho s(N)}{4} \right\}$$

$$+ \mathbb{P}\left\{ \max_{1 \leq m \leq t_{N(\beta)}} \left\| \sum_{j=1}^{m} \mathbb{G}_j(f) \right\|_{\mathcal{F}} > \frac{\rho s(N)}{4} \right\}$$

$$+ \sum_{k=N(\beta)}^{N-1} \mathbb{P}\left\{ \max_{t_k+1 \leq m \leq t_{k+1}} \left\| \sum_{j=t_k+1}^{m} [f(X_j) - \mathbb{E}f(X)] \right\|_{\mathcal{F}} > \frac{\rho s(N)}{8} \right\}$$

$$+ \sum_{k=N(\beta)}^{N-1} \mathbb{P}\left\{ \max_{t_k+1 \leq m \leq t_{k+1}} \left\| \sum_{j=t_k+1}^{m} \mathbb{G}_j(f) \right\|_{\mathcal{F}} > \frac{\rho s(N)}{8} \right\}$$

$$+ \mathbb{P}\left\{ \max_{N(\beta) \leq j < N} \left\| \sum_{k=N(\beta)}^{j} \left( \sqrt{n_k} \alpha_{n_k}^{(k)} - \sqrt{n_k} \mathbb{G}^{(k)} \right) \right\|_{\mathcal{F}} > \frac{\rho s(N)}{4} \right\} =: \sum_{i=1}^{5} P_i(\rho, N).$$

It is easy to show using inequalities (5.9) and (5.10), along with the choice of $1/2 < \beta = \alpha/(1+\alpha) < 1$, that for any $\gamma > 0$ for all large enough $\rho$,

$$(5.14) \qquad \sum_{i=1}^{2} P_i(\rho, N) \leq t_N^{-\gamma}/4, \text{ for all } N \geq 1.$$

For instance, consider $P_1(\rho, N)$. Observe that

$$P_1(\rho, N) \leq \mathbb{P}\left\{ \max_{1 \leq m \leq t_{N(\beta)}} \sqrt{m} \|\alpha_m\|_{\mathcal{F}} > C\sqrt{t_{N(\beta)}}(B + \tau_N) \right\},$$

where

$$\tau_N = \left( \frac{\rho s(N)}{4} - B \right) / \left( C\sqrt{t_{N(\beta)}} \right).$$

Now $\sqrt{t_{N(\beta)}} \sim c_2 N^{\alpha/2}$ for some $c_2 > 0$. Therefore by (5.13) for some $c_3 > 0$,

$$\tau_N \sim c_3 N^{1-\tau_1 \alpha} (\log N)^{\tau_2}.$$

Since by (5.11) we have $1 - \tau_1 \alpha > 0$, we readily get from inequality (5.9) that for any $\gamma > 0$ and all large enough $\rho$, $P_1(\rho, N) \leq t_N^{-\gamma}/8$, for all $N \geq 1$. In the same way we get using inequality (5.10) that for any $\gamma > 0$ and all large enough $\rho$, $P_2(\rho, N) \leq t_N^{-\gamma}/8$, for all $N \geq 1$. Hence we have (5.14).

In a similar fashion one can verify that for any $\gamma > 0$ and all large enough $\rho$,

$$(5.15) \qquad \sum_{i=3}^{4} P_i(\rho, N) \leq t_N^{-\gamma}/4, \text{ for all } N \geq 1.$$

To see this, notice that

$$P_3(\rho, N) \leq N\mathbb{P}\left\{ \max_{1 \leq m \leq n_N} \sqrt{m} \|\alpha_m\|_{\mathcal{F}} > \rho s(N)/8 \right\}$$

and

$$P_4(\rho, N) \leq N\mathbb{P}\left\{ \max_{1 \leq m \leq n_N} \left\| \sum_{j=1}^{m} \mathbb{G}_j(f) \right\|_{\mathcal{F}} > \rho s(N)/8 \right\}.$$

Since $\sqrt{n_N} \sim N^{\alpha/2}$ and $N \sim c_3 t_N^{1/(\alpha+1)}$ for some $c_3 > 0$, we get (5.15) by proceeding as above using inequalities (5.9) and (5.10).

Next, recalling the definition of $s(N)$, we get

$$P_5(\rho, N) \leq \mathbb{P}\left\{ \sum_{k=N(\beta)}^{N} \left\| \sqrt{n_k} \alpha_{n_k}^{(k)} - \sqrt{n_k} \mathbb{G}^{(k)} \right\|_{\mathcal{F}} > \frac{\rho s(N)}{4} \right\}$$

$$\leq \sum_{k=N(\beta)}^{N} \mathbb{P}\left\{ \left\| \sqrt{n_k} \alpha_{n_k}^{(k)} - \sqrt{n_k} \mathbb{G}^{(k)} \right\|_{\mathcal{F}} > \frac{\rho n_k^{1/2-\tau_1} (\log n_k)^{\tau_2}}{4} \right\},$$

which by (5.12) for any $\lambda > 0$ and $\rho = \rho(\alpha, \lambda) > 0$ large enough is

$$\leq N \left( \left[ N^\beta \right]^\alpha \right)^{-\lambda}, \text{ for all } N \geq 1,$$

which, in turn, for large enough $\lambda > 0$ is $\leq t_N^{-\gamma}/2$. Thus for all $\gamma > 0$ there exists a $\rho > 0$ so that

$$\sum_{i=1}^{5} P_i(\rho, N) \leq t_N^{-\gamma}, \text{ for all } N \geq 1.$$

Since $\alpha$ can be any number satisfying $1/2 < \tau_1 \alpha < 1$ and $t_{N+1}/t_N \to 1$, this implies (2.4) for $\rho = \rho(\alpha, \lambda)$ large enough. The almost sure statement (2.5) follows trivially from (2.4) using a simple blocking and the Borel–Cantelli lemma on the just constructed probability space. This proves Theorem 1. □

*Proof of Theorem 2.* The proof follows along the same lines as that of Theorem 1. Therefore for the sake of brevity we shall only outline the proof. Here we borrow ideas from the proof of Theorem 6.2 of Dudley and Philipp [14]. Recall that in Theorem 2 we assume that $1/2 < r_0 < 1$ in Proposition 2, which means that $0 < \kappa := (1 - r_0)/2r_0 < 1/2$. For $k \geq 1$ set

(5.16) $$t_k = \left[ \exp\left( k^{1-\kappa} \right) \right] \text{ and } n_k = t_k - t_{k-1}, \text{ where } t_0 = 1.$$

Now for some $b > 0$ we get $n_k \sim b^2 k^{-\kappa} t_k$,

$$\frac{\sqrt{n_k}}{(\log n_k)^\kappa} \sim \frac{b\sqrt{t_k}}{k^{\kappa(1-\kappa)+\kappa/2}} = \frac{b\sqrt{t_k}}{k^{\kappa+\theta}},$$

where $\theta = \kappa \left( \frac{1}{2} - \kappa \right) > 0$. Choose $0 < \beta < 1$ and set $N(\beta) = \left[ N^\beta \right]$. Using an integral approximation we get for suitable constants $c_1 > 0$ and $c_2 > 0$, for all large $N$

(5.17) $$\frac{c_1 \sqrt{t_N}}{N^\theta} \leq s(N) := \sum_{k=N(\beta)}^{N} \frac{\sqrt{n_k}}{(\log n_k)^\kappa} \leq \frac{c_2 \sqrt{t_N}}{N^\theta} \leq \frac{c_2 \sqrt{t_N}}{(\log(t_N))^{\theta/(1-\kappa)}}.$$

Also for all large $N$,

(5.18) $$s(N)/\sqrt{n_N} \geq \frac{c_1}{2b} N^{\kappa/2 - \kappa\left(\frac{1}{2} - \kappa\right)} =: c_0 N^{\kappa^2}.$$

For later use note that for any $0 < \beta < 1$ and $\zeta > 0$

(5.19) $$\frac{s(N)}{\sqrt{t_{N(\beta)}} N^\zeta} \to \infty, \text{ as } N \to \infty,$$

and observe that

(5.20) $$t_{N+1}/t_N \to 1, \text{ as } N \to \infty.$$

Constructing a probability space and defining $P_i(\rho, N)$, $i = 1, \ldots, 5$, as in the proof of Theorem 1, but with $n_k$, $t_k$ and $s(N)$ as given in (5.16) and (5.17) the proof now goes much like that of Theorem 1. In particular, using inequalities (5.9) and (5.10), and noting that $N \sim (\log(t_N))^{1/(1-\kappa)}$, one can check that for some $\nu > 0$, for all large enough $N$,

$$\sum_{i=1}^{4} P_i(\rho, N) \leq \exp\left( -(\log(t_N))^\nu \right)$$

and by arguing as in the proof of Theorem 1, but now using Proposition 2, we easily see that for every $H > 0$ there is a probability space on which sit i.i.d. $X_1, X_2...$, and i.i.d. $\mathbb{G}_1, \mathbb{G}_2, \ldots$, and a $\rho > 0$ such that

$$P_5(\rho, N) \leq (\log(t_N))^{-H-1}, \text{ for all } N \geq 1.$$

Since for all $H > 0$,

$$\log(t_N)^H \left( \exp\left(-(\log(t_N))^\nu\right) + (\log(t_N))^{-H-1} \right) \to 0, \text{ as } N \to \infty,$$

this in combination with (5.17) and (5.20) proves that (2.8) holds with $\tau = \theta/(1-\kappa)$ and $\rho(\tau, H)$ large enough. A simple blocking argument shows that (2.9) follows from (2.8). Choose $H > 1$ in (2.8). Notice that for any $k \geq 1$,

$$\mathbb{P}\left\{ \bigcup_{2^k < n \leq 2^{k+1}} \left\{ \max_{1 \leq m \leq n} \left\| \sqrt{m}\alpha_m - \sum_{i=1}^m \mathbb{G}_i \right\|_{\mathcal{F}} > \sqrt{2n}\rho(\tau, H)(\log n)^{-\tau} \right\} \right\}$$
$$\leq \mathbb{P}\left\{ \max_{1 \leq m \leq 2^{k+1}} \left\| \sqrt{m}\alpha_m - \sum_{i=1}^m \mathbb{G}_i \right\|_{\mathcal{F}} > \sqrt{2^{k+1}}\rho(\tau, H)(\log 2^{k+1})^{-\tau} \right\}$$
$$\leq ((k+1)\log 2)^{-H}.$$

Hence (2.9) holds by the Borel-Cantelli lemma. $\square$

**Acknowledgement.** The authors thank the referee for pointing out a number of oversights and misprints, as well as showing the way to an improvement in Theorem 2.